\documentstyle[11pt,leqno]{article}
\newtheorem{theorem}{{\sc Theorem}}
\newcommand{\bt}{\begin{theorem}}
\newcommand{\et}{\end{theorem}}
\newcommand{\newsection}[1]{\setcounter{equation}{0} \setcounter{theorem}{0}
\section{#1}}

\newcommand{\NI}{\noindent}
\newcommand{\bea}{\begin{eqnarray}}
\newcommand{\eea}{\end{eqnarray}}

\def \spec#1 {\mathop{#1}}

\def \b #1 {\bf #1}

\newcommand {\CC}{\centerline}

\newcommand{\clf}{{\cal F}}

\newcommand{\ity}{\infty}
\newcommand{\raro}{\rightarrow}

\newcommand{\vsp}{\vskip 1em}
\newcommand{\vspp}{\vskip 2em}

\newcommand{\be}{\begin{equation}}
\newcommand{\ee}{\end{equation}}
\newcommand{\ben}{\begin{eqnarray*}}
\newcommand{\een}{\end{eqnarray*}}

\begin{document}
\sloppy
\CC {\Large{\bf  Nonparametric Estimation  of Linear Multiplier in }}
\CC {\Large {\bf Stochastic  Differential Equations Driven by $\alpha$-Stable Noise }}
\vsp 
\CC {\bf B.L.S. Prakasa Rao}
\CC{\bf CR Rao Advanced Institute of Mathematics, Statistics }
\CC{\bf and Computer Science, Hyderabad, India}
\vspp 
\NI{\bf Abstract:} We discuss nonparametric estimation of linear multiplier in a trend coefficient in models
governed by a stochastic differential equation driven by an $\alpha$-stable small noise. 
\vsp 
\NI{\bf Keywords and phrases}: Stochastic differential equation; Trend coefficient; Linear multiplier; Nonparametric estimation; Kernel method; Levy process. 
\vsp 
\NI AMS Subject classification (2010): Primary 62M09; Secondary 60G52.

\newsection{Introduction}

Let $(\Omega, \clf, P)$ be a probability space equipped with a right continuous and increasing family of $\sigma$-algebras $\{\clf_t, t \geq 0\}.$ Let $\{Z_t, t \geq 0\}$ be a standard $\alpha$-stable Levy process with $Z_1$ distributed as $S_\alpha(1,\beta,0).$ A random variable $Z$ is said to have a stable distribution $S_\alpha(\sigma, \beta,\mu)$ with index of stability $\alpha \in (0,2],$ scale parameter $\sigma \in (0,\ity),$ skewness parameter $\beta \in [-1.1],$ and location parameter $\mu \in (-\ity,\ity)$ if it has the characteristic function (c.f) of the following form:
\bea
\phi_Z(u) &=& E[\exp(iuZ)]\\\nonumber
&=& \exp \{-\sigma^\alpha |u|^\alpha(1-i\beta \;sgn(u) \tan (\frac{\alpha \pi}{2}))+i \mu u\} \;\;\mbox{if}\;\; \alpha \neq 1\\\nonumber
&=& \exp \{ -\sigma|u|(1+i\beta \frac{2}{\pi} \;sgn(u) \log |u|)+i\mu u\} \;\;\mbox{if}\;\; \alpha = 1.\\\nonumber
\eea
If $\mu=0,$ then we say that the random variable $Z$ is strictly stable. If, in addition $\beta=0,$ we say that the random variable $Z$ is symmetric $\alpha$-stable (cf. Samorodnitsky and Taqqu (1994), Sato (1999)). Here after we assume that $1<\alpha<2.$ Suppose that $X=\{X_t, 0\leq t \leq T\}$ is a stochastic process satisfying the stochastic differential equation (SDE)
\be
dX_t= \theta(t) X_t dt+ \epsilon \;dZ_t, 0\leq t \leq T, X_0=x_0
\ee
where $Z=\{Z_t, 0\leq t \leq T\}$ is a standard $\alpha$-stable Levy process with $Z_1$ distributed as $S_\alpha(1,\beta,0).$ Suppose that $\alpha$ and $\beta$ are known with $1 <\alpha<2$ but the linear multiplier $\theta(.)$ is unknown. The problem is to estimate the function 
$\theta(.)$ and the trend function based on the observations of the process $X$ over the interval $[0,T].$ Parameter estimation for Ornstein-Uhlenbeck process driven by $\alpha$-stable Levy motion is investigated in Hu and Long (2007, 2009). Least squares estimation for discretely observed Ornstein-Uhlenbeck processes with small Levy noises is studied in Long (2009) and Ma (2010). Parametric estimation for a class of SDE driven by small stable noises from discrete observations is studied in Long (2010). Shen and Xu (2014) studied estimation of the drift parameter for SDEs driven by small noises which is more than pure jump $\alpha$-stable noises. Shen et al. (2018) discussed parameter estimation for Ornstein-Uhlenbeck processes driven by a fractional Levy process. Nonparametric estimation of the trend function for stochastic differential equations driven by a mixed fractional Brownian motion or a sub-fractional Brownian motion are investigated in Prakasa Rao (2019, 2020).\\
\vsp 
\newsection{Levy processes}

An ${\clf_t}$-adapted stochastic process $\{Z_t, t \geq 0\}$ ts called $\alpha$-stable Levy process if (i) $Z_0=0$ a.s.,(ii) $Z_t-Z_s$ has the distribution $S_\alpha((t-s)^{1/\alpha}, \beta,0)$ for $0\leq s < t < \ity,$ and (iii) for any $0 \leq t_0 < t_1<\dots<t_m < \ity,$ the random variables $Z_{t_0},Z_{t_1}-Z_{t_0}, \dots, Z_{t_m}-Z_{t_{m-1}}$ are independent. Ito-type stochastic integrals with respect to $\alpha$-stable Levy processes were investigated in Kallenberg (1975, 1992) and Rosinski and Woyczynski (1986). Let $L^{\alpha}$ be the family of all real-valued ${\clf_t}$-predictable processes on $\Omega \times [0,\ity)$ such that for every $T>0, \int_0^T|\phi(t,\omega)|^\alpha dt < \ity$ a.s. It is known that a predictable process $\phi$ is integrable with respect to a strictly $\alpha$-stable Levy process $Z$, that is, $\int_0^T\phi(t)dZ_t$ exists for every $T>0,$ if and only if $\phi \in L^\alpha.$ The following moment inequality  is due to Long (2010) improving the result in Theorem 3.2 in Rosinski and Woyczynski (1985).\\

\NI{\bf Theorem 2.1:} {\it Suppose $Z$ is a $\alpha$-stable Levy process and $\phi$ is a predictable process such that $\phi \in L^\alpha.$ Then there exists a positive constant $c$ depending only on $\alpha$ and $\beta$ independent of $T$ such that}
\be
E[\sup_{0\leq t \leq T}|\int_0^t\phi(s)dZ_s|] \leq c E[(\int_0^T|\phi(t)|^\alpha dt)^{1/\alpha}].
\ee

In particular, by choosing $\phi(s)\equiv 1,$ it follows that
\be
E[\sup_{0\leq t \leq T}|Z_t|] \leq c T^{1/\alpha}
\ee
where $c$ is a positive constant depending on $\alpha.$

The following result, due to Rosinski and Woyczynski (1985) and Kallenberg (1992, Theorem 4.1), gives a method for finding the distribution of a stochastic integral with respect to an $\alpha $-stable Levy process. Let $\phi_{+}$ and $\phi_{-}$ denote the positive and negative parts of a function $\phi.$

\NI{\bf Theorem 2.2:} {\it Suppose $\phi \in L^{\alpha}$ and $Z$ is a strictly $\alpha$-stable Levy process. Then (i) there exists independent processes $Z^\prime$ and $Z^{\prime \prime}$ with the same distribution as $Z$ such that
\be
\int_0^t\phi(s)dZ_s= Z^\prime(\int_0^t [\phi_{+}(s)]^\alpha ds) -Z^{\prime\prime}(\int_0^t[\phi_{-}(s)]^\alpha ds)
\ee
almost surely. (ii) If $Z$ is symmetric, that is $\beta=0,$ then there exists a $\alpha$-stable Levy process $Z^\prime$ such that $Z^\prime$ and $Z$ have the same distribution and
\be
\int_0^t\phi(s)dZ_s= Z^\prime(\int_0^t|\phi(s)|^\alpha ds) 
\ee
almost surely.}
\vsp
\newsection{Preliminaries}

Let us consider the stochastic differential equation
\be
dX_t= \theta(t) X_t\; dt + \epsilon \;  dZ_t, X_0=x_0 , 0 \leq t \leq T
\ee
where the trend function $\theta(.)$ is unknown. We would like to estimate the function $\theta(.)$ and the drift $\theta(t)X_t$  based on the
observation $X= \{ X_t, 0 \leq t \leq T\}. $ Let $x= \{x_t, 0 \leq t \leq T\}$ be the solution of the ordinary differential equation
\be
\frac{dx_t}{dt}=\theta(t) x_t , x_0 , 0 \leq t \leq T.
\ee
Observe that
$$x_t= x_0 \exp(\int_0^t\theta(s)ds).$$
\vsp
We assume that (A1) the function $\theta(t)$ is bounded over the interval $[0,T]$ by a constant $L$. 
\vsp
\NI{\bf Lemma 3.1: } {\it Let  $ X_t$ and $x_t$ be  the solutions of the equation (3.1) and
(3.2) respectively. Then, with probability one,}
\be
(a)|X_t-x_t|\leq  e^{Lt} \epsilon \sup_{0\leq s \leq t}|Z_s|
\ee
and
\be
(b)\sup_{0 \leq t \leq T} E|X_t-x_t| \leq e^{LT} \epsilon T^{1/\alpha}.
\ee
\vsp
\NI{\bf Proof of (a) :} Let $u_t=|X_t-x_t| $. Then
\bea
u_t & \leq & \int^t_0 \left| \theta(v)(X_v-x_v) \right| dv + \epsilon \;|Z_t|\\\nonumber
& \leq & L \int^t_0 u_v dv + \epsilon \;\sup_{0\leq s \leq t}|Z_s|.
\eea
Applying the Gronwall's lemma (cf. Lemma 1.12, Kutoyants (1994), p.26),  it follows that
\be
u_t \leq  \epsilon \sup_{0\leq s \leq t}|Z_s| e^{Lt}. 
\ee
\vsp
\NI{\bf Proof of (b) :} From (3.3), we have ,
\bea
E|X_t-x_t| & \leq & e^{Lt} \epsilon \; E (|\sup_{0\leq s \leq t}|Z_s|) \\\nonumber
& = & e^{Lt}\epsilon t^{1/\alpha}\\\nonumber
\eea
by the moment inequality in Theorem 2.1 applied to the function $\phi \equiv 1.$ Hence
\be
\sup_{0 \leq t \leq T} E |X_t-x_t|  \leq e^{LT}\epsilon T^{1/\alpha}. \\
\ee
\vsp
\newsection{Estimation of the Drift function}

Let $\Theta_0(L)$ denote the class of all functions $\theta(.) $ with the same bound $L$. Let $\Theta_k(L) $ denote the class of all functions $\theta(.) $ which are uniformly bounded by the same constant $L$ and which are $k$-times differentiable satisfying the condition
$$|\theta^{(k)}(x)-\theta^{(k)}(y)|\leq L^\prime |x-y|, x,y \in R$$
for some constant $L^\prime >0.$ Here $g^{(k)}(x)$ denotes the $k$-th derivative of $g(.)$ at $x$ for $k \geq 0.$ If $k=0,$ we interpret the function $g^{(0)}$ as the function  $g(.).$\\
Let $G(u)$ be a  bounded function  with finite support $[A,B]$ with $A<0<B$ satisfying the condition 

\NI{$(A_2)$}$ G(u) =0\;\; \mbox{for}\;\; u <A \;\;\mbox{and} u>B;\;\;\mbox{and} \int^B_A G(u) du =1.$

Boundedness of the function $G(.)$ with finite support $[A,B]$ implies that

$$\int_-\ity^\ity |G(u)|^\alpha du <\ity; \int_\ity^\ity |u^jG(u)| du <\ity, j \geq 0.$$

We define a kernel type estimator $\widehat{\theta}_t$ of the function $\theta(t)$ by the relation 
\be
\widehat{\theta}_t X_t = \frac{1}{\varphi_\epsilon}\int^T_0 G \left(\frac{\tau-t}{\varphi_\epsilon} \right) d X_\tau
\ee
where the normalizing function  $ \varphi_\epsilon \rightarrow 0 $ as  $ \epsilon \rightarrow 0. $ Let $E_\theta(.)$ denote the expectation when the function $\theta(.)$ is the linear multiplier.
\vsp
\NI{\bf Theorem 4.1:}  {\it Suppose that the linear multiplier $\theta(.) \in \Theta_0(L)$ and  the function
\ $ \varphi_\epsilon \rightarrow 0$  such that $ \epsilon \varphi^{-1}_\epsilon \longrightarrow  0 $  as $\epsilon
\rightarrow  0$. Further suppose that the conditions $(A_1), (A_2)$   hold. Then, for any $ 0 < c
\leq d < T ,$ the estimator $\widehat{\theta}_t$ satisfies the property}
\be
\lim_{\epsilon \rightarrow 0} \sup_{\theta(.) \in \Theta_0(L)} \sup_{c\leq t \leq d } E_\theta ( |\widehat{\theta}_t X_t - \theta(t)x_t|)= 0.
\ee
\vsp
In addition to the conditions $(A_1),(A_2),$ suppose the following condition holds.\\
\NI{$(A_3)$}$ \int^\infty_{-\infty} u^j G(u)   du = 0 \;\;\mbox{for}\;\; j=1,2,...k.$
\vsp
\NI{\bf Theorem 4.2:} {\it Suppose that the function $ \theta(.) \in \Theta_{k+1}(L)$ and  $
\varphi_\epsilon = \epsilon^{\frac{1}{k-\frac{1}{\alpha}+2}}.$ Suppose the conditions $(A_1)-(A_3)$ hold. Then}
\be
\limsup_{\epsilon \rightarrow 0} \sup_{\theta(.) \in \Theta_{k+1}(L)}\sup_{c \leq t \leq d} E_\theta (| \widehat{\theta}_t X_t- \theta(t)x_t|)\epsilon^{\frac{-(k+1)}{k-\frac{1}{\alpha}+2}}  < \infty.
\ee
\vsp
\NI{\bf Theorem 4.3:} {\it Suppose that the function $\theta(.) \in \Theta_{k+1}(L)$ and $ \varphi_\epsilon= \epsilon^{\gamma}$ where $\gamma= \frac{1}{(k+2)-\frac{1}{\alpha}}.$ Suppose the conditions $(A_1)-(A_3)$ hold. Then the asymptotic distribution of
$$ \varphi_\epsilon^{-(k+1)}(\widehat{\theta}_t X_t- \theta(t)x_t), $$
as $\epsilon \raro 0$ is the distribution of the random variable
$$(\int_0^T G(u)_+^{\alpha}du)^{1/\alpha}U_1-(\int_0^TG(u)_-^{\alpha}du)^{1/\alpha}U_2$$
where $U_1$ and $U_2$ are independent random variables with $\alpha $- stable distribution $S_\alpha(1,\beta,0)$ shifted by the constant $m$ defined by
$$m= \frac{J^{(k+1)}(x_t)}{(k+1)!}\int_{-\ity}^{\ity}G(u)u^{k+1}du$$
and $J(t)=\theta(t)x_t.$}
\vsp
\newsection{Proofs of Theorems}

\NI{\bf Proof of Theorem 4.1 :} From the equation (3.1), we have
\bea
\;\;\;\\\nonumber
E_\theta[|\widehat{\theta}_tX_t -\theta(t)x_t|]  &=& E_\theta \{ |\frac{1}{\varphi_\epsilon}  \int^T_0 G \left(\frac{\tau-t}{\varphi_\epsilon} \right) \left(\theta(\tau)X_\tau -\theta(\tau) x_\tau \right)  d \tau \\ \nonumber
& &+ \frac{1}{\varphi_\epsilon}
 \int^T_0 G \left(\frac{\tau-t}{\varphi_\epsilon}\right) \theta(\tau) x_\tau d \tau- \theta(t) x_t
 + \frac{\epsilon}{\varphi_\epsilon} \int^T_0 G \left(\frac{\tau-t}{\varphi_\epsilon} \right)
 dZ_\tau|\}\\ \nonumber
 & \leq  &  E_\theta \left[ |\frac{1}{\varphi_\epsilon}  \int^T_0 G \left(\frac{\tau-t}{\varphi_\epsilon} \right) (\theta(\tau)X_\tau -\theta(\tau)x_\tau) d\tau |\right] \\ \nonumber
 & & + E_\theta \left[|\frac{1}{\varphi_\epsilon} \int^T_0 G \left(\frac{\tau-t}{\varphi_\epsilon} \right)\theta(\tau) x_\tau d\tau -\theta(t) x_t |\right] \\ \nonumber
 & & +  \frac{\epsilon}{\varphi_{\epsilon}} E_\theta \left[ |\int^T_0 G \left(  \frac{\tau-t}{\varphi_\epsilon}\right) d Z_\tau|\right]\\ \nonumber
 &= & I_1+I_2+I_3 \;\;\mbox{(say).}\;\;\\ \nonumber
\eea
Apply the change of variables $u= (t-\tau)\varphi_{\epsilon}^{-1}$ and denote $\epsilon_1= \min(\epsilon^\prime, \epsilon^{\prime\prime})$, where $\epsilon^\prime= \sup\{\epsilon: \varphi_\epsilon \leq -\frac{c}{A}\}$ and $\epsilon^{ \prime\prime}= \sup\{\epsilon: \varphi_\epsilon \leq -\frac{T-d}{B}\}.$ Then, for $\epsilon <\epsilon_1,$
\bea
\;\;\;\\ \nonumber
I_3 &= & \frac{ \epsilon}{\varphi_\epsilon} E_\theta \left( \int_0^T G\left(\frac{\tau-t}{\varphi_\epsilon} \right) d Z_\tau|
\right) \\ \nonumber
& \leq & \frac{ \epsilon}{\varphi_\epsilon} C_1
 \left[\int^T_0 \left\{ |G \left(\frac{\tau-t}{\varphi_\epsilon}
\right)|\right\}^{\alpha} d \tau  \right]^{1/\alpha} \\\nonumber
&& \;\;\;\;\; \mbox{( cf . Long (2010) )} \\ \nonumber
&=& \frac{ \epsilon}{\varphi_\epsilon}[\int_{-\ity}^{\ity} |G(u)|^\alpha du]^{1/\alpha}( \varphi_\epsilon)^{1/\alpha}\\\nonumber
& \leq & \frac{C_2 \epsilon}{\varphi_\epsilon} [\varphi^{1/\alpha}_\epsilon ](\;\; \mbox{(by using $(A_2)$  )} \\ \nonumber
\eea
for some positive constant $C_2$ depending on $T,\alpha, L, A$ and $B.$ Since $\epsilon \varphi_\epsilon^{-1} \raro 0$ and $\varphi_\epsilon \raro 0$ as $\epsilon \raro 0$, it follows that $I_1$ tends to zero as $\epsilon \raro 0.$
\vsp
Furthermore
\bea
\;\;\;\\ \nonumber
I_2 &= & E_\theta \left[ |\frac{1}{\varphi_\epsilon} \int^T_0 G\left(
\frac{\tau-t}{\varphi_\epsilon}\right) \theta(\tau) x_\tau) d \tau - \theta(t) x_t|
\right] \\ \nonumber
& \leq &  E_\theta \left[ \int^\infty_{-\infty} |G(u)
\left(\theta( t+\varphi_\epsilon u) x_{t+\varphi_\epsilon u}-\theta(t) x_t \right)  |\ du \right]
\\ \nonumber
& \leq &  L\left[ \int^\infty_{-\infty} |G(u) u| \varphi_\epsilon du
\right]\\ \nonumber
& \leq  & C_3 \varphi_\epsilon \\ \nonumber
\eea
for some positive constant $C_3$ depending on $T,\alpha, L,L_1, A $ and $B$.  Hence $I_2$ tends to zero as $\epsilon \raro 0.$ Furthermore note that
\bea
\;\;\;\\\nonumber
I_1 &= & E_\theta \left[ |\frac{1}{\varphi_\epsilon} \int^T_0 G
\left(\frac{\tau-t}{\varphi_\epsilon} \right) (\theta (\tau)X_\tau -\theta(\tau) x_\tau)
d\tau |\right] \\\nonumber
&= & E_\theta  \left[ |\int^\infty_{-\infty} G(u)
\left( \theta(t+\varphi_\epsilon u) X_{t+\varphi_\epsilon u}  - \theta(t+\varphi_\epsilon u) x_{t+\varphi_\epsilon u}
)\right) du|\right]\\\nonumber
& \leq &  (\int^\infty_{-\infty} |G(u)| L E_\theta \left(
|X_{t+\varphi_\epsilon u}-x_{t+\varphi_\epsilon u} |\right) \ du
\;\;\mbox{(by using the condition $(A_1)$)}\\\nonumber
& \leq &  \int^\infty_{-\infty} |G(u)| \;\;L \sup_{0 \leq t +
\varphi_\epsilon u \leq T}E_\theta \left(|X_{t+\varphi_\epsilon u}
-x_{t+\varphi_\epsilon u}|\right) \ du \\\nonumber
& \leq & C_4 \epsilon \;\;\mbox{(by using (3.4))}\\\nonumber
\eea
for some positive constant $C_4$ depending on $T, \alpha$ and $L$. Hence $I_3$  tends to zero as  $\epsilon \raro 0.$  Theorem 4.1 is now proved by using the equations (5.1) to (5.4).
\vsp
\NI{\bf Remarks:} From the proof presented above, it is possible to choose the functions $c_\epsilon $ and 
$d_\epsilon $ such that $c_\epsilon \raro 0, d_\epsilon \raro T$ 
and satisfy the conditions 
$$\frac{c_\epsilon}{\varphi_\epsilon}\geq -A, \frac{T-d_\epsilon}{\varphi_\epsilon}\geq B$$
(for instance, choose $c_\epsilon = -A \varphi_\epsilon  $ and $d_\epsilon = T-B \varphi_\epsilon $. Then the estimator $\hat \theta_t$ satisfies the property that
\be
\lim_{\epsilon \rightarrow 0} \sup_{\theta(t) \in \Theta_0(L)} \sup_{c_\epsilon \leq t \leq d_\epsilon } E_\theta ( |\widehat{\theta}_tX_t - \theta(t)x_t|)= 0.
\ee
\vsp
\NI {\bf Proof of Theorem 4.2 :}
Let $ J(t) = \theta(t) x_t.$ By the Taylor's formula, for any $u \in R,$

$$ J(y) = J(u) +\sum^k_{r=1} J^{(r)}(u) \frac{(y-u)^r}{r!} +[ J^{(k)}(z)-J^{(k)}(u)] \frac{(y-u)^k}{k!} $$
for some $z$ such that $|z-u|\leq |y-u|.$ Using this expansion, the equation (3.2) and the conditions  in the expression $I_2$ defined in the proof of  Theorem 4.1, it follows that
\ben
\;\;\\\nonumber
I_2 & \leq &  \left[|\int^\infty_{-\infty} G(u) \left(J(t+\varphi_\epsilon u) - J(t) \right)  \ du |\right]\\\nonumber
&= & [ |\sum^k_{j=1} J^{(j)}(t) (\int^\infty_{-\infty}G(u) u^j du )\varphi^j_\epsilon (j!)^{-1}\\\nonumber
& & \;\;\;\;+(\int^\infty_{-\infty}G(u) u^k (J^{(k)}(z_u) -J^{(k)}(x_t))du \varphi^k_\epsilon (k !)^{-1}|]\\ \nonumber
\een
for some $z_u$ such that $|x_t-z_u|\leq |x_{t+\varphi_\epsilon u}-x_t| \leq C|\varphi_\epsilon u|$ for some positive constant $C.$ Hence
\bea
I_2 & \leq & C_5 L \left[  \int^\infty_{-\infty} |G(u)u^{k+1}|\varphi^{k+1}_\epsilon (k!) ^{-1}  du  \right]\\\nonumber
& \leq & C_6 (k!)^{-1} \varphi^{k+1}_\epsilon \int^\infty_{-\infty} |G(u) u^{k+1}| du  \\\nonumber
&\leq & C_7 \varphi_\epsilon^{k+1}\\ \nonumber
\eea
for some positive constant $C_7$ depending on $A, B,  \alpha,T$ and $L$. Combining the relations (5.2), (5.4) and (5.6) , we get that there exists a positive constant $C_8$ depending on $\alpha,T,L,A, B$ such that 
$$ \sup_{c \leq t \leq d}E_\theta|\widehat{\theta}_tx_t-\theta(t) x_t| \leq C_{8} (\epsilon\varphi^{-1}_\epsilon \varphi^{1/\alpha}_\epsilon+  \varphi^{k+1}_\epsilon +\epsilon). $$ 
Choosing $ \varphi_\epsilon = \epsilon^{\frac{1}{k+2-\frac{1}{\alpha}}},$  we get that 
$$ \limsup_{\epsilon \rightarrow 0} \sup_{\theta(.) \in \Theta_{k+1}(L) } \sup_{c \leq t \leq d} E_\theta |\widetilde{\theta}_t X_t - \theta(t) x_t|\epsilon^ {-\frac{(k+1)}{{k-\frac{1}{\alpha}+2}}} < \infty. $$ 
This completes the proof of Theorem 4.2. 
\vsp 
\vsp \NI{\bf Proof of Theorem 4.3:}

From the equation (3.1), we obtain that
\bea
\lefteqn{ \widehat{\theta}_tX_t -\theta(t) x_t}\\\nonumber
 &= &[ \frac{1}{\varphi_\epsilon} \int^T_0 G \left(\frac{\tau-t}{\varphi_\epsilon} \right)
 \left( \theta(\tau) X_\tau- \theta(\tau) x_\tau \right) \  d \tau \\\nonumber
 & & + \frac{1}{\varphi_\epsilon} \int^T_0 G \left( \frac{\tau-t}{\varphi_\epsilon}\right) \theta(\tau) x_\tau d\tau -\theta(t)x_t+ \frac{\epsilon}{\varphi_\epsilon} \int^T_0 G \left( \frac{\tau-t}{\varphi_\epsilon}\right) dZ_\tau]\\\nonumber
 &= & [ \int^\infty_{-\infty} G(u) (\theta(t+\varphi_\epsilon u) X_{t+\varphi_\epsilon u} - \theta (t+\varphi_\epsilon u) x_{t+\varphi_\epsilon u}) \ du  \\\nonumber
 & & +\int^\infty_{-\infty} G(u) (\theta (t+\varphi_\epsilon u) x_{t+\varphi_\epsilon u}- \theta(t) x_t) \ du \\\nonumber
 & &+ \frac{\epsilon}{\varphi_{\epsilon}}\int^T_0 G\left(\frac{\tau-t}{\varphi_\epsilon}
 \right) d Z_\tau].\\\nonumber
\eea
Let $J(t)= \theta(t)x_t.$ By the Taylor's formula, for any $u \in R,$
$$ J(y) = J(u) +\sum^{k+1}_{r=1} J^{(r)}(u) \frac{(y-u)^j}{j!} +[J^{(k+1)}(z)-J^{(k+1)}(x)] \frac{(y-u)^{k+1}}{(k+1)!} $$
for some $z$ such that $|z-u|\leq |y-u|.$ Let 
$$m= \frac{J^{(k+1)}(x_t)}{(k+1)!}\int_{-\ity}^{\ity}G(u)u^{k+1}du$$
and
$$ R_1(t)= \varphi^{-(k+1)}_\epsilon\int^T_0 G\left(\frac{\tau-t}{\varphi_\epsilon} \right) (\theta(\tau) X_\tau -\theta (\tau) x_\tau)
d\tau.$$
By arguments similar to those given in (5.2) for obtaining upper bounds, it follows that 
$$E|R_1(t)| \leq C \varphi_\epsilon^{-k}\epsilon.$$
Let 
$$ R_2(t)= \varphi_\epsilon^{-(k+1)} \int^T_0 G\left(\frac{\tau-t}{\varphi_\epsilon} \right) (\theta(\tau) x_\tau -\theta (t) x_t)
d\tau.$$
Observe that
$$R_2(t)=m+o(1)$$
by an application of the Taylor's expansion under the condition $(A_3).$ Furthermore
\bea
\varphi_\epsilon^{-(k+1)}(\hat \theta_t x_t-\theta(t) x_t) &= & \epsilon \varphi^{-(k+2)}_\epsilon\int^T_0 G\left(\frac{\tau-t}{\varphi_\epsilon} \right) d Z_\tau + R_2(t) + R_1(t)\\\nonumber
&=&  \epsilon \varphi^{-(k+2)}_\epsilon\int^T_0 G\left(\frac{\tau-t}{\varphi_\epsilon} \right) d Z_\tau + m+o(1)+O_p(\varphi_\epsilon^{-k}\epsilon).\\\nonumber 
\eea
Let $\varphi_\epsilon$ be chosen so that $(\varphi_\epsilon)^{-1/\alpha}= \epsilon \varphi_\epsilon^{-(k+2)}.$ One such choice is $\varphi_\epsilon =\epsilon^v$ where $v= (k+2-\frac{1}{\alpha})^{-1}.$ 
We will now study the asymptotic behaviour of the random variable 
$$W_\epsilon= (\varphi_\epsilon)^{-1/\alpha}\int^T_0 G\left(\frac{\tau-t}{\varphi_\epsilon} \right) d Z_\tau $$
as $\epsilon \raro 0.$  Note that	
\be
\varphi_\epsilon^{-(k+1)}(\hat \theta_t X_t-\theta(t) x_t)=W_\epsilon+ m+o_p(1).
\ee
By Theorem 2.2, there exist two independent processes $Z^\prime$ and $Z^{\prime\prime}$ with the same distribution as $Z$ such that
\bea
\lefteqn{\int_0^T G\left(\frac{\tau-t}{\varphi_\epsilon} \right) d Z_\tau}\\\nonumber
&= & Z^\prime(\int_0^T(G\left(\frac{\tau-t}{\varphi_\epsilon} \right))_+^{\alpha}d\tau)-Z^{\prime\prime}(\int_0^T(G\left(\frac{\tau-t}{\varphi_\epsilon} \right))_-^{\alpha}d\tau)\\\nonumber
&=& (\int_0^T(G\left(\frac{\tau-t}{\varphi_\epsilon} \right))_+^{\alpha}d\tau)^{1/\alpha}U_1-(\int_0^T(G\left(\frac{\tau-t}{\varphi_\epsilon} \right))_-^{\alpha}d\tau)^{1/\alpha}U_2\\\nonumber
&= &\phi_\epsilon^{1/\alpha} (\int_0^TG(u)_+^{\alpha}du)^{1/\alpha}U_1-\phi_\epsilon^{1/\alpha} (\int_0^TG(u)_-^{\alpha}du)^{1/\alpha}U_2
\eea
where $U_1$ and $U_2$ are two independent random variables with $\alpha $- stable distribution $S_\alpha(1,\beta,0).$ Hence the asymptotic distribution of the random variable
 $$\varphi_\epsilon^{-(k+1)}(\hat \theta_t X_t-\theta(t) x_t)$$
is the distribution of the random variable
$$(\int_0^T G(u)_+^{\alpha}du)^{1/\alpha}U_1-(\int_0^TG(u)_-^{\alpha}du)^{1/\alpha}U_2$$
where $U_1$ and $U_2$ are independent random variables with $\alpha $- stable distribution $S_\alpha(1,\beta,0)$ shifted by the constant $m$ defined by
$$m= \frac{J^{(k+1)}(x_t)}{(k+1)!}\int_{-\ity}^{\ity}G(u)u^{k+1}du$$
where $J(t)=\theta(t)x_t.$
\vsp
Note that the results given above deal with asymptotic  properties of the estimator for the function 
$$J(t)= \theta(t)x_t= \theta(t) x_0 \exp(\int_0^t \theta(s)\;ds).$$ 
\vsp
We will now present another method for the estimation of the linear multiplier $\theta(t).$
\newsection{Estimation of the Multiplier $\theta(.)$}

Let $\Theta_\rho(L_\gamma)$ be a class of functions uniformly bounded  and $k$-times continuously differentiable for some integer 
$k \geq 1 $ with the $k$-th derivative satisfying the Holder condition of the order $\gamma \in (0,1):$
$$|\theta^{(k)}(t)-\theta^{(k)}(s)|\leq L_\gamma |t-s|^\gamma, \rho=k+\gamma.$$

From the Lemma 3.1, it follows that
$$|X_t-x_t| \leq \epsilon e^{Lt}\sup_{0\leq s \leq T} |Z_s|.$$
Let
$$A_t= \{\omega: \inf_{0\leq s \leq t}X_s(\omega)\geq \frac{1}{2}x_0e^{-Lt}\}$$
and let $A=A_T.$ Define the process $Y$ with the differential
$$dY_t=\theta(t) I(A_t) dt + \epsilon X_t^{-1}I(A_t)\;dZ_t, 0\leq t \leq T.$$
We will now construct an estimator of the function $\theta(.)$ based on the observation of the process $Y$ over the interval $[0,T].$
Define the estimator
$$\tilde \theta(t)= I(A) \frac{1}{\varphi_\epsilon}\int_0^T G(\frac{t-s}{\varphi_\epsilon})dY_s$$
where the kernel function $G(.)$ satisfies the conditions $(A1)-(A3)$. Observe that
\ben
E|\tilde \theta(t)-\theta(t)| &= & E|I(A) \frac{1}{\varphi_\epsilon}\int_0^T G(\frac{t-s}{\varphi_\epsilon})(\theta(s)-\theta(t))ds\\\nonumber
&& \;\;\;\; + I(A^c)\theta(t)+I(A)\frac{\epsilon}{\varphi_\epsilon}\int_0^T G(\frac{t-s}{\varphi_\epsilon})X_s^{-1}dZ_s|\\\nonumber
&\leq & E|I(A)\int_R G(u)[\theta(t+u\varphi_\epsilon)-\theta(t)]du|+ |\theta(t)| P(A^c)\\\nonumber
&& \;\;\;\; + \frac{\epsilon}{\varphi_\epsilon}|E[I(A)\int_0^T G(\frac{t-s}{\varphi_\epsilon})X_s^{-1}dZ_s]|\\\nonumber
&=& I_1+I_2+I_3. \;\;\mbox{(say)}.\\\nonumber
\een
Applying the Taylor's theorem and using the fact that the function $\theta(t)\in \Theta_\rho(L_\gamma)$, it follows that
\ben
I_1 \leq \frac{L_\gamma}{(k+1)!}\varphi_\epsilon^\rho \int_R|G(u)u^\rho|du.
\een
Note that, by Lemma 3.1,
\ben
P(A^c) & = & P(\inf_{0\leq t \leq T}X_t < \frac{1}{2}x_0e^{-LT})\\\nonumber
&\leq & P(\inf_{0 \leq t \leq T}|X_t-x_t| + \inf_{0\leq t \leq T}x_t < \frac{1}{2}x_0e^{-LT})\\\nonumber
&\leq & P(\inf_{0 \leq t \leq T}|X_t-x_t| < -\frac{1}{2}x_0e^{-LT})\\\nonumber
&\leq & P(\sup_{0 \leq t \leq T}|X_t-x_t| > \frac{1}{2}x_0e^{-LT})\\\nonumber
&\leq & P(\epsilon e^{LT}\sup_{0 \leq t \leq T}|Z_t|>\frac{1}{2}x_0e^{-LT})\\\nonumber
&= & P(\sup_{0 \leq t \leq T}|Z_t|>\frac{x_0}{2\epsilon}e^{-2LT})\\\nonumber
&\leq & \frac{D}{\alpha (2-\alpha)}T(\frac{x_0}{2\epsilon}e^{-2LT})^{-\alpha}\\\nonumber
\een
for some positive constant $D$ by a maximal inequality for stable stochastic integrals in Gine and Marcus (1983) (cf. Joulin (2006)). The upper bound obtained above and the fact that $|\theta(s)|\leq L, 0\leq s \leq T$ leads  an upper bound for the term $I_2.$ We have used the inequality
$$x_t= x_0 \exp(\int_0^t\theta (s)ds)\geq x_0 e^{-Lt}$$
in the computations given above. Applying Theorem 2.1, it follows that
\ben
\lefteqn{|E[I(A)\int_0^T G(\frac{t-s}{\varphi_\epsilon})X_s^{-1}dZ_s|}\\\nonumber
&=& |E[\int_0^T G(\frac{t-s}{\varphi_\epsilon})X_s^{-1}I(A_s)dZ_s|\\\nonumber
&\leq & C E[(\int_0^T|G(\frac{t-s}{\varphi_\epsilon})X_s^{-1}I(A_s)|^\alpha ds)^{1/\alpha}]\\\nonumber
&\leq &C e^{LT} (\int_0^T|G(\frac{t-s}{\varphi_\epsilon})|^\alpha ds)^{1/\alpha}\\\nonumber
&=& Ce^{LT} (\varphi_\epsilon)^{1/\alpha} [\int_R|G(u)|^{\alpha}du]^{1/\alpha}\\\nonumber
\een
for some positive constant $C$ depending on $\alpha$ which leads to an upper bound on the term $I_3.$ Combining the above estimates, it follows that
\ben
E|\tilde \theta(t)-\theta(t)|\leq C_1\varphi_\epsilon^\rho + C_2 \epsilon^{\alpha}+ C_3\epsilon \varphi_\epsilon^{(1/\alpha)-1}
\een
for some positive constants $C_i, i=1,2,3$ depending on $\alpha$ and $L.$ Choosing $\varphi_\epsilon=\epsilon^{\alpha/\rho},$ we obtain that
\ben
E|\tilde \theta(t)-\theta(t)|\leq C_4 \epsilon^{(\rho-\alpha+1)/\rho}+ C_5 \epsilon^{\alpha}
\een
for some positive constants $C_4$ and $C_5.$ It easy to see that $0< \frac{\rho-\alpha+1}{\rho} <\alpha$ for all $\rho> \alpha-1$ which holds  since $1<\alpha<2$ and $\rho >1.$ Hence we obtain the following result implying the uniform consistency of the estimator $\tilde \theta(t)$ as an estimator of $\theta(t)$ as $\epsilon \raro 0.$\\

\NI{\bf Theorem 6.1:} {\it Let $\theta \in \Theta_\rho(L)$ and $\varphi_\epsilon= \epsilon^{\alpha/(\rho}.$ Suppose the conditions $(A_1)-(A_3)$ hold. Then, for any interval $[c,d] \subset [0,T],$ }
\ben
\limsup_{\epsilon \raro 0}\sup_{\theta(.)\in \Theta_\rho(L)}\sup_{c\leq t \leq d}E|\tilde \theta(t)-\theta(t)| \epsilon^{-\frac{\rho+1-\alpha}{\rho}}<\ity.
\een
\vsp
\NI{\bf Funding :} This work was supported under the scheme ``INSA Senior Scientist" by the Indian National Science Academy  while the author was at the CR Rao Advanced Institute for Mathematics, Statistics and Computer Science, Hyderabad 500046, India.\\
\vsp
\NI{\bf References :}
\begin{description}
\item Gine, E. and Marcus, M. (1983) The central limit theorem for stochastic integrals with respect to Levy processes, {\it Ann. Probability},  {\bf 11}, 58-77.

\item Hu, Y. and Long, H. (2007) Parameter estimation for Ornstein-Uhlenbeck processes driven by $\alpha$-stable Levy motions, {\it Communications on Stochastic Analysis}, {\bf 1}, 175-192.

\item Hu, Y. and Long, H. (2009) Least squares estimator for Ornstein-Uhlenbeck processes driven by $\alpha$-stable Levy motions, {\it Stochastic Process Appl.}, {\bf 119}, 2465-2480.

\item Joulin, A. (2006) On maximal inequalities for stable stochastic integrals, Prepint, Department of Mathematics, Universite de La Rochelle.

\item Kallenberg, O. (1975) On the existence and path properties of stochastic integrals, {\it Ann. Probab.}, {\bf 3}, 262-280.

\item Kallenberg, O. (1992) Some time change representations of stable integrals via predictable transformations of local martingale, {\it Stochastic Process appl.}, {\bf 40}, 199-223.

\item Kutoyants, Y.A.(1994) {\it Identification of Dynamical Systems with Small Noise}. Kluwer: Dordrecht.

\item Long, H. (2009) Least squares estimator for discretely observed  Ornstein-Uhlenbeck processes with small Levy noises, {\it Statist. Probab. Lett.}, {\bf 79}, 2076-2085.

\item Long, H. (2010) Parameter estimation for a class of stochastic differential equations driven by small stable noises from discrete observations, {\it Acta Mathematica Scientia}, {\bf 30B}, 645-663.

\item Ma, C. (2010) A note on ``Least squares estimator for discretely observed  Ornstein-Uhlenbeck processes with small Levy noises", {\it Statist. Probab. Lett.}, {\bf 80}, 1528-1531.

\item  Prakasa Rao, B.L.S. (2019) Nonparametric estimation of trend for stochastic differential equations driven by mixed fractional Brownian motion, {\it Stoch. Anal. and Appl.}, {\bf 37}, 271-280.

\item  Prakasa Rao, B.L.S. (2020) Nonparametric estimation of trend for stochastic differential equations driven by sub-fractional Brownian motion, {\it Random Oper. Stoch. Equ.} (to appear).

\item Rosinski, J. and Woyczynski, W.A. (1985) Moment inequalities for real and vector $p$-stable  stochastic integrals, {\it Probability Theory in Banach Spaces V}, Lecture Notes in Math. Vol. 1153: 369-386, Springer.

\item Rosinski, J. and Woyczynski, W.A. (1986) On Ito stochastic integration with respect to $p$-stable motion: inner clock, integrability of sample paths, double and multiple integrals, {\it Ann. Probab.}, {\bf 14}, 271-286.

\item  Samorodnitsky, G. and Taqqu, M.S. (1994) {\it Stable non-Gaussian Random Processes: stochastic Models with Infinite Variance}, Chapman and Hall, New York.

\item Sato, K.I. ( 1999) {\it Levy Processes and Infinitely Divisible Distributions}, Cambridge University Press, Cambridge. 

\item Shen, L. and Xu, Q. (2014) Statistical inference for stochastic differential equations with small noises, {\it Abstract and Applied analysis}, Volume 2014, Article ID 473681, 6 pages.

\item Shen, G., Li, Y. and Gao, Z. (2018) Parameter estimation for Ornstein-Uhlenbeck process driven by fractional Levy process, {\it Journal of Inequalities and Applications}, Volume 2018, 2018: 356.

\end{description}
\end{document}